\documentclass[12pt]{article}

\usepackage{amsmath, epsfig, cite}
\usepackage{amssymb}
\usepackage{amsfonts}
\usepackage{latexsym}
\usepackage{float}

\newtheorem{thm}{Theorem}[section]

\newtheorem{cor}[thm]{Corollary}

\newtheorem{lem}[thm]{Lemma}

\numberwithin{equation}{section}

\newcommand{\qed}{{\hfill$\square$}\medskip}

\newcommand*\pFqskip{8mu}
\catcode`,\active
\newcommand*\pFq{\begingroup
        \catcode`\,\active
        \def ,{\mskip\pFqskip\relax}%
        \dopFq
}
\catcode`\,12
\def\dopFq#1#2#3#4#5{%
        {}_{#1}F_{#2}\biggl[\genfrac..{0pt}{}{#3}{#4};#5\biggr]%
        \endgroup
}

\begin{document}


\begin{center}
{\Large\bf Some supercongruences on truncated ${}_3F_2$\\[5pt]
hypergeometric series }
\end{center}

\vskip 2mm \centerline{Ji-Cai Liu}
\begin{center}
{\footnotesize Department of Mathematics, Wenzhou University, Wenzhou 325035, PR China\\
{\tt jcliu2016@gmail.com} }
\end{center}


\vskip 0.7cm \noindent{\bf Abstract.}
In 2003, Rodriguez-Villegas conjectured four supercongruences on the truncated ${}_3F_2$ hypergeometric series for certain modular K3 surfaces, which were gradually proved by several authors. Motivated by some supercongruences on combinatorial numbers such as Ap\'ery numbers and Domb numbers, we establish some new supercongruences on the truncated ${}_3F_2$ hypergeometric series, which extend the four Rodriguez-Villegas supercongruences.

\vskip 2mm
\noindent{\it MR Subject Classifications}: Primary 11A07, 11S80; Secondary 33C20, 33B15

\section{Introduction}
In 2003, Rodriguez-Villegas \cite{rv-b-2003} observed the remarkable relationship between the number of points over $F_p$ on certain Calabi-Yau manifolds and truncated hypergeometric series.
In doing so, he conjectured numerically 22 supercongruences for hypergeometric Calabi-Yau manifolds of dimension $D\le 3$.

To state these results, we first define the truncated hypergeometric series. For complex numbers $a_i, b_j$ and $z$, with none of the $b_j$ being negative integers or zero, the truncated hypergeometric series are defined as
\begin{align*}
\pFq{r}{s}{a_1,a_2,\cdots,a_r}{b_1,b_2,\cdots,b_s}{z}_n=\sum_{k=0}^{n-1}\frac{(a_1)_k (a_2)_k\cdots (a_r)_k}{(b_1)_k (b_2)_k\cdots (b_s)_k}\cdot \frac{z^k}{k!},
\end{align*}
where $(a)_0=1$ and $(a)_k=a(a+1)\cdots (a+k-1)$ for $k\ge 1$.

Throughout this paper let $p\ge5$ be a prime.
For manifolds of dimension $D=2$, Rodriguez-Villegas \cite{rv-b-2003} conjectured four supercongruences, which are related to certain modular K3 surfaces. These are all of the form
\begin{align}
\pFq{3}{2}{\frac{1}{2},-a,a+1}{1,1}{1}_{p}\equiv c_p\pmod{p^2},\label{na-1}
\end{align}
where $a=-\frac{1}{2}, -\frac{1}{3},-\frac{1}{4},-\frac{1}{6}$ and $c_p$ is the $p$-th Fourier coefficient of a weight three modular form on a congruence subgroup of $SL(2, \mathbb{Z})$.
The case when $a=-\frac{1}{2}$ was confirmed by van Hamme \cite{vanhamme-b-1987}, Ishikawa \cite{ishikawa-nmj-1990} and Ahlgren \cite{ahlgren-b-2001}.
The other cases when $a=-\frac{1}{3},-\frac{1}{4},-\frac{1}{6}$ were partially proved by
Mortenson \cite{mortenson-pams-2005}, and finally proved by Sun \cite{sunzw-aa-2012}.

Recall that the $p$-adic Gamma function \cite[Definition 11.6.5]{cohen-b-2007} is defined as
\begin{align*}
\Gamma_p(x)=\lim_{m\to x}(-1)^m\prod_{\substack{1< k < m\\
(k,p)=1}}k,
\end{align*}
where the limit is for $m$ tending to $x$ $p$-adically in $\mathbb{Z}_{\ge 0}$.
Let $\langle a\rangle_p$ denote the least non-negative integer $r$ with $a\equiv r \pmod{p}$.

Sun \cite[Theorem 2.5]{sunzh-jnt-2014} and the author \cite[Theorem 1.2]{liujc-a07686-2016}
proved that for any $p$-adic integer $a$,
\begin{align}
&\hskip -5mm\pFq{3}{2}{\frac{1}{2},-a,a+1}{1,1}{1}_{p}\notag\\
&\equiv
\begin{cases}
0,\quad&\text{if $\langle a\rangle_p\equiv 1\pmod{2}$}\\[5pt]
(-1)^{\frac{p+1}{2}}\Gamma_p\left(-\frac{a}{2}\right)^2\Gamma_p\left(\frac{a+1}{2}\right)^2 ,\quad&\text{if $\langle a\rangle_p\equiv 0\pmod{2}$}
\end{cases}\pmod{p^2},\label{na-2}
\end{align}
which extends \eqref{na-1} for any $p$-adic integer $a$.
Guo and Zeng \cite[Theorem 1.3]{gz-aa-2015} have established an interesting
$q$-analogue of \eqref{na-2} when $\langle a\rangle_p\equiv 1\pmod{2}$.

In the ingenious proof of the irrationality of $\zeta(3)$, Ap\'ery  \cite{apery-asterisque-1979} introduced the numbers
$$A_n=\sum_{k=0}^n{n\choose k}^2{n+k\choose k}^2,$$
which are known as Ap\'ery numbers.  A motivating example in this paper is the following supercongruence:
\begin{align*}
A_{np}\equiv A_n\pmod{p^3},
\end{align*}
which was originally conjectured by Chowla et al. \cite{ccc-jnt-1980}, and first proved by
Gessel \cite{gessel-jnt-1982}.
Some similar types of supercongruences on combinatorial numbers such as Almkvist-Zudilin numbers, Domb numbers and Ap\'ery-like numbers have been studied by several authors, see for example, Amdeberhan and Tauraso \cite{at-aa-2016}, Chan, Cooper and Sica \cite{ccs-ijnt-2010}, Osburn and Sahu \cite{os-aam-2011}, and Osburn, Sahu and Straub \cite{oss-pems-2016}.

In this paper, we aim to establish the same type of supercongruences on truncated ${}_3F_2$ hypergeometric series. The first aim of this paper is to prove the following result, which strengthens \cite[Theorem 2.5]{sunzh-jnt-2014}.

\begin{thm}\label{nt-2}
Let $p\ge 5$ be a prime and $n$ be a positive integer. For any $p$-adic integer $a$ with $\langle a\rangle_p\equiv 1 \pmod{2}$, we have
\begin{align}
\pFq{3}{2}{\frac{1}{2},-a,a+1}{1,1}{1}_{np}\equiv 0 \pmod{p^2}.\label{na-4}
\end{align}
\end{thm}

The second aim of this paper is to extend the supercongruence \eqref{na-1} for $\langle a\rangle_p\equiv 0 \pmod{2}$.
\begin{thm}\label{nt-3}
Let $p\ge 5$ be a prime and $n$ be a positive integer. For $a\in\left\{-\frac{1}{2},-\frac{1}{3},-\frac{1}{4},-\frac{1}{6}\right\}$ with $\langle a\rangle_p\equiv 0 \pmod{2}$, we have
\begin{align}
&\hskip -5mm\pFq{3}{2}{\frac{1}{2},-a,a+1}{1,1}{1}_{np}\notag\\
&\equiv (-1)^{\frac{p+1}{2}}\Gamma_p\left(-\frac{a}{2}\right)^2\Gamma_p\left(\frac{a+1}{2}\right)^2
\pFq{3}{2}{\frac{1}{2},-a,a+1}{1,1}{1}_{n} \pmod{p^2}.\label{na-5}
\end{align}
\end{thm}

We see that letting $n=1$ in \eqref{na-5} reduces to \eqref{na-1} for $\langle a\rangle_p\equiv 0 \pmod{2}$.

{\noindent \bf Remark.} Numerical calculation via {\tt Sigma} suggests that supercongruence \eqref{na-5} cannot be extended to any $p$-adic integer $a$ with $\langle a\rangle_p\equiv 0 \pmod{2}$ in the direction of \eqref{na-2}.

Replacing $n$ by $p^{r-1}$ in \eqref{na-5} and then using induction, we get the following result.
\begin{cor}
Let $p\ge 5$ be a prime and $r$ be a positive integer. For $a\in\left\{-\frac{1}{2},-\frac{1}{3},-\frac{1}{4},-\frac{1}{6}\right\}$ with $\langle a\rangle_p\equiv 0 \pmod{2}$, we have
\begin{align*}
\pFq{3}{2}{\frac{1}{2},-a,a+1}{1,1}{1}_{p^r}
\equiv (-1)^{\frac{r(p+1)}{2}}\Gamma_p\left(-\frac{a}{2}\right)^{2r}\Gamma_p\left(\frac{a+1}{2}\right)^{2r}
\pmod{p^2}.
\end{align*}
\end{cor}

The proof of Theorem \ref{nt-2} and \ref{nt-3} will be given in Section 2 and 3, respectively.

\section{Proof of Theorem \ref{nt-2}}
In order to prove Theorem \ref{nt-2}, we first establish the following two lemmas.
\begin{lem}
For any odd integer $n$, we have
\begin{align}
&\sum_{k=0}^n{2k\choose k}^2{n+k\choose 2k}\left(-\frac{1}{4}\right)^k=0,\label{nb-1}\\
&\sum_{k=0}^n{2k\choose k}^2{n+k\choose 2k}\left(-\frac{1}{4}\right)^k(2H_{2k}-3H_k)=0,\label{nb-2}\\
&\sum_{k=0}^n{2k\choose k}^2{n+k\choose 2k}\left(-\frac{1}{4}\right)^k(2H_{n+k}-H_k)=0,\label{nb-3}
\end{align}
where $H_n=\sum_{k=1}^n\frac{1}{k}$ denotes the $n$-th harmonic number.
\end{lem}
{\noindent \it Proof.}
An identity \cite[(1), p.16]{bailey-b-1964} says
\begin{align}
&\pFq{3}{2}{a,b,c}{\frac{1}{2}(a+b+1),2c}{1}\notag\\
&=\frac{\Gamma\left(\frac{1}{2}\right)\Gamma\left(\frac{1}{2}+c\right)\Gamma\left(\frac{1}{2}+\frac{1}{2}a+\frac{1}{2}b\right)
\Gamma\left(\frac{1}{2}-\frac{1}{2}a-\frac{1}{2}b+c\right)}{\Gamma\left(\frac{1}{2}+\frac{1}{2}a\right)
\Gamma\left(\frac{1}{2}+\frac{1}{2}b\right)\Gamma\left(\frac{1}{2}-\frac{1}{2}a+c\right)\Gamma\left(\frac{1}{2}-\frac{1}{2}b+c\right)}.
\label{nb-4}
\end{align}

Letting $a=-n,b=n+1$ and $c=\frac{1}{2}$ in \eqref{nb-4} reduces to
\begin{align*}
\pFq{3}{2}{-n,n+1,\frac{1}{2}}{1,1}{1}=\frac{\Gamma\left(\frac{1}{2}\right)^2}
{\Gamma\left(\frac{1-n}{2}\right)^2\Gamma\left(\frac{n+2}{2}\right)^2}.
\end{align*}
Since $(1-n)/2$ is always a non-positive integer, we have $\displaystyle\lim_{x\to \frac{1-n}{2}}1/\Gamma\left(x\right)^2=0$, and so
\begin{align*}
\pFq{3}{2}{-n,n+1,\frac{1}{2}}{1,1}{1}=0.
\end{align*}
Noting that
\begin{align}
\frac{\left(\frac{1}{2}\right)_k}{(1)_k}&=\frac{{2k\choose k}}{4^k},\label{nb-5}\\
\frac{(-n)_k(n+1)_k}{(1)_k^2}&=(-1)^k{2k\choose k}{n+k\choose 2k},\label{nb-6}
\end{align}
we are led to \eqref{nb-1}.

Let $A_n$ and $B_n$ respectively denote the numbers
\begin{align*}
A_n=\sum_{k=0}^{2n-1}{2k\choose k}^2{2n-1+k\choose 2k}\left(-\frac{1}{4}\right)^k(2H_{2k}-3H_k),
\end{align*}
and
\begin{align*}
B_n=\sum_{k=0}^{2n-1}{2k\choose k}^2{2n-1+k\choose 2k}\left(-\frac{1}{4}\right)^k(2H_{2n-1+k}-H_k).
\end{align*}
Using the software package {\tt Sigma} developed by Schneider \cite{schneider-slc-1999}, we find that $A_n$ and $B_n$ respectively satisfy the following recurrences:
\begin{align*}
&16n^3(n+1)(4n+5)A_n-4(n+1)(2n+1)(4n+3)(4n^2+6n+1)A_{n+1}\\
&+(2n+1)(2n+3)^3(4n+1)A_{n+2}=0,
\end{align*}
and
\begin{align*}
16n^3(n+1)B_n-4(n+1)(2n+1)^3B_{n+1}+(2n+1)^2(2n+3)^2B_{n+2}=0.
\end{align*}
It is easy to verify that $A_n=B_n=0$ for $n=1,2$. This implies that
$A_n=B_n=0$ for all $n\ge 1$, and so \eqref{nb-2} and \eqref{nb-3} clearly hold.
\qed

\begin{lem}
Suppose $p\ge 5$ is a prime.
For any $p$-adic integer $a$ with $\langle a\rangle_p\equiv 1 \pmod{2}$, we have
\begin{align}
&\sum_{k=0}^{p-1}{2k\choose k}^2{a+k\choose 2k}\left(-\frac{1}{4}\right)^k\sum_{i=1}^{k}\frac{1}{-a-1+i}\notag\\
&\equiv \sum_{k=0}^{p-1}{2k\choose k}^2{a+k\choose 2k}\left(-\frac{1}{4}\right)^k\sum_{i=1}^{k}\frac{1}{a+i}\notag\\
&\equiv \frac{1}{2}\sum_{k=0}^{p-1}{2k\choose k}^2{a+k\choose 2k}\left(-\frac{1}{4}\right)^kH_k\pmod{p}.\label{nb-7}
\end{align}
\end{lem}
{\noindent \it Proof.}
From \eqref{nb-1} and \eqref{nb-3}, we obtain that for any odd integer $n$,
\begin{align*}
\sum_{k=0}^n{2k\choose k}^2{n+k\choose 2k}\left(-\frac{1}{4}\right)^k(2H_{n+k}-2H_n-H_k)=0,
\end{align*}
which is
\begin{align}
\sum_{k=0}^n{2k\choose k}^2{n+k\choose 2k}\left(-\frac{1}{4}\right)^k\left(\sum_{i=1}^k\frac{1}{n+i}-\frac{1}{2}H_k\right)=0.\label{nb-8}
\end{align}
Since $\langle a\rangle_p$ is an odd integer, it follows from \eqref{nb-8} that
\begin{align*}
&\sum_{k=0}^{p-1}{2k\choose k}^2{a+k\choose 2k}\left(-\frac{1}{4}\right)^k\sum_{i=1}^{k}\frac{1}{a+i}\\
&\equiv \sum_{k=0}^{\langle a\rangle_p}{2k\choose k}^2{\langle a\rangle_p+k\choose 2k}\left(-\frac{1}{4}\right)^k\sum_{i=1}^{k}\frac{1}{\langle a\rangle_p+i}\pmod{p}\\
&=\frac{1}{2}\sum_{k=0}^{\langle a\rangle_p}{2k\choose k}^2{\langle a\rangle_p+k\choose 2k}\left(-\frac{1}{4}\right)^kH_k \\
&\equiv \frac{1}{2}\sum_{k=0}^{p-1}{2k\choose k}^2{a+k\choose 2k}\left(-\frac{1}{4}\right)^kH_k\pmod{p}.
\end{align*}

Let $b=p-\langle a\rangle_p$. It is clear that
$a\equiv -b \pmod{p}$ and $b-1$ is an odd integer with $0\le b-1\le p-1$. It follows that
\begin{align}
&\sum_{k=0}^{p-1}{2k\choose k}^2{a+k\choose 2k}\left(-\frac{1}{4}\right)^k\sum_{i=1}^{k}\frac{1}{-a-1+i}\notag\\
&\equiv\sum_{k=0}^{p-1}{2k\choose k}^2{-b+k\choose 2k}\left(-\frac{1}{4}\right)^k\sum_{i=1}^{k}\frac{1}{b-1+i}\pmod{p}.\label{nb-9}
\end{align}
Since ${-b+k\choose 2k}={b-1+k\choose 2k}$, by \eqref{nb-8} and \eqref{nb-9} we get
\begin{align*}
\text{LHS \eqref{nb-9}}&\equiv\sum_{k=0}^{b-1}{2k\choose k}^2{b-1+k\choose 2k}\left(-\frac{1}{4}\right)^k\sum_{i=1}^{k}\frac{1}{b-1+i}\pmod{p}\\
&=\frac{1}{2}\sum_{k=0}^{b-1}{2k\choose k}^2{b-1+k\choose 2k}\left(-\frac{1}{4}\right)^kH_k\\
&\equiv \frac{1}{2}\sum_{k=0}^{p-1}{2k\choose k}^2{-a-1+k\choose 2k}\left(-\frac{1}{4}\right)^kH_k\pmod{p}\\
&=\frac{1}{2}\sum_{k=0}^{p-1}{2k\choose k}^2{a+k\choose 2k}\left(-\frac{1}{4}\right)^kH_k .
\end{align*}
This completes the proof of \eqref{nb-7}.
\qed

{\noindent\it Proof of Theorem \ref{nt-2}.}
By \eqref{nb-5} and \eqref{nb-6}, we can rewrite \eqref{na-4} as
\begin{align*}
\sum_{k=0}^{np-1}{2k\choose k}^2{a+k\choose 2k}\left(-\frac{1}{4}\right)^k\equiv 0 \pmod{p^2}.
\end{align*}
We first prove that for any non-negative integer $r$,
\begin{align}
\sum_{k=rp}^{(r+1)p-1}{2k\choose k}^2{a+k\choose 2k}\left(-\frac{1}{4}\right)^k\equiv 0 \pmod{p^2},\label{nb-10}
\end{align}
which is equivalent to
\begin{align}
\sum_{k=0}^{p-1}{2k+2rp\choose k+rp}^2{a+k+rp\choose 2k+2rp}\left(-\frac{1}{4}\right)^{k+rp}\equiv 0 \pmod{p^2}.\label{nb-11}
\end{align}

Applying the following congruence \cite[(2.1)]{liujc-jdea-2016}:
\begin{align}
{2rp+2k\choose rp+k}\equiv {2r\choose r}{2k\choose k}(1+2rp(H_{2k}-H_k))\pmod{p^2},\label{nx}
\end{align}
we immediately get
\begin{align}
{2rp+2k\choose rp+k}^2\equiv {2r\choose r}^2{2k\choose k}^2(1+4rp(H_{2k}-H_k))\pmod{p^2}.\label{nb-12}
\end{align}
By \eqref{nb-12} and the fact that ${2k\choose k}^2\equiv 0\pmod{p^2}$ for $\frac{p-1}{2}< k \le p-1$, we have
\begin{align}
&\text{LHS \eqref{nb-11}}\notag\\
&\equiv \sum_{k=0}^{\frac{p-1}{2}}{2r\choose r}^2{2k\choose k}^2{a+k+rp\choose 2k+2rp}\left(-\frac{1}{4}\right)^{k+rp}(1+4rp(H_{2k}-H_k)) \pmod{p^2}.\label{nb-13}
\end{align}

Let $\delta$ denote the number $\delta=(a-\langle a\rangle_p)/p$. It is clear that $\delta$ is a $p$-adic integer and $a=\langle a\rangle_p+\delta p$.
Note that
\begin{align*}
&{a+k+rp\choose 2k+2rp}\notag\\
&={\langle a\rangle_p+k+(\delta+r) p\choose 2k+2rp}\notag\\
&={\langle a\rangle_p+(\delta +r)p\choose 2rp}\prod_{i=1}^k(\langle a\rangle_p+(\delta +r)p+i)\prod_{i=1}^k(\langle a\rangle_p+(\delta -r)p+1-i)\prod_{i=1}^{2k}(2rp+i)^{-1}.
\end{align*}
It is easy to see that for $0\le k\le\frac{p-1}{2}$,
\begin{align*}
\prod_{i=1}^{2k}(2rp+i)^{-1}\equiv \prod_{i=1}^{2k}\left(\frac{1}{i}-\frac{2rp}{i^2}\right)\equiv \frac{1-2rpH_{2k}}{(2k)!}\pmod{p^2},
\end{align*}
and
\begin{align*}
&\prod_{i=1}^k(\langle a\rangle_p+(\delta +r)p+i)\prod_{i=1}^k(\langle a\rangle_p+(\delta -r)p+1-i)\\
&\equiv\prod_{i=1}^{2k}\left(\langle a\rangle_p-k+i\right)\cdot\left(1+p\left(\sum_{i=1}^k\frac{\delta+r}{\langle a\rangle_p+i}+\sum_{i=1}^k\frac{\delta-r }{\langle a\rangle_p+1-i}\right)\right)\pmod{p^2}.
\end{align*}
It follows that for $0\le k \le \frac{p-1}{2}$,
\begin{align}
{a+k+rp\choose 2k+2rp}
&\equiv {\langle a\rangle_p+(\delta +r)p\choose 2rp}
{\langle a\rangle_p+k\choose 2k}\notag\\
&\times\left(1+p\left(\sum_{i=1}^k\frac{r+\delta}{\langle a\rangle_p+i}+\sum_{i=1}^k\frac{r-\delta }{-\langle a\rangle_p-1+i}-2rH_{2k}\right)\right)\pmod{p^2}.\label{nb-14}
\end{align}
Substituting \eqref{nb-14} into \eqref{nb-13} gives
\begin{align}
\text{LHS \eqref{nb-11}}
&\equiv \left(-\frac{1}{4}\right)^{rp}{2r\choose r}^2{\langle a\rangle_p+(\delta +r)p\choose 2rp}\sum_{k=0}^{\frac{p-1}{2}}{2k\choose k}^2{\langle a\rangle_p+k\choose 2k}\left(-\frac{1}{4}\right)^k\notag\\
&\times \left(1+p\left(\sum_{i=1}^k\frac{r+\delta}{\langle a\rangle_p+i}+\sum_{i=1}^k\frac{r-\delta }{-\langle a\rangle_p-1+i}+2rH_{2k}-4rH_k\right)\right)\pmod{p^2}.\label{nb-15}
\end{align}

Note that ${2k\choose k}^2\equiv 0\pmod{p^2}$ for $\frac{p-1}{2}< k \le p-1$ and $\langle a\rangle_p$ is an odd integer. It follows from \eqref{nb-1} that
\begin{align}
\sum_{k=0}^{\frac{p-1}{2}}{2k\choose k}^2{\langle a\rangle_p+k\choose 2k}\left(-\frac{1}{4}\right)^k
\equiv \sum_{k=0}^{\langle a\rangle_p}{2k\choose k}^2{\langle a\rangle_p+k\choose 2k}\left(-\frac{1}{4}\right)^k
=0\pmod{p^2}.\label{nb-16}
\end{align}
On the other hand, using \eqref{nb-7} and \eqref{nb-2} we have
\begin{align}
&\sum_{k=0}^{\frac{p-1}{2}}{2k\choose k}^2{\langle a\rangle_p+k\choose 2k}\left(-\frac{1}{4}\right)^k\left(\sum_{i=1}^k\frac{r+\delta}{\langle a\rangle_p+i}+\sum_{i=1}^k\frac{r-\delta }{-\langle a\rangle_p-1+i}+2rH_{2k}-4rH_k\right)\notag\\
&\equiv \sum_{k=0}^{\langle a\rangle_p}{2k\choose k}^2{\langle a\rangle_p+k\choose 2k}\left(-\frac{1}{4}\right)^k\left(2rH_{2k}-3rH_k\right)
=0\pmod{p}.\label{nb-17}
\end{align}
Then the proof of \eqref{nb-11} directly follows from \eqref{nb-15}-\eqref{nb-17}.

Taking the sum over $r$ from $0$ to $n-1$ on both sides of \eqref{nb-10} gives
\begin{align*}
\sum_{k=0}^{np-1}{2k\choose k}^2{a+k\choose 2k}\left(-\frac{1}{4}\right)^k\equiv 0 \pmod{p^2},
\end{align*}
which is equivalent to \eqref{na-4}.
This completes the proof of Theorem \ref{nt-2}. \qed

\section{Proof of Theorem \ref{nt-3}}
The Fermat quotient of an integer $a$ with respect to an odd prime $p$ is given by
\begin{align*}
q_p(a)=\frac{a^{p-1}-1}{p}.
\end{align*}
\begin{lem} \label{l1} (Eisenstein)
Suppose $p$ is an odd prime and $r$ is a positive integer. For non-zero $p$-adic integers $a$ and $b$, we have
\begin{align*}
q_p(ab)&\equiv q_p(a)+q_p(b)\pmod{p},\\
q_p(a^r)&\equiv rq_p(a) \pmod{p}.
\end{align*}
\end{lem}

\begin{lem}\label{l2} (Lehmer \cite{lehmer-am-1938})
For any prime $p\ge 5$, we have
\begin{align*}
H_{\lfloor p/2\rfloor}\equiv -2q_p(2)\pmod{p},
\end{align*}
where $\lfloor x \rfloor$ denotes the greatest integer less than or equal to a real number $x$.
\end{lem}

\begin{lem} (See \cite[Lemma 2.4]{liujc-bams-2017})
Let $p\ge 5$ be a prime, $r$ and $k$ be non-negative integers. For $0\le k \le p-1$ and $a\in\left \{-\frac{1}{2},-\frac{1}{3},-\frac{1}{4},-\frac{1}{6}\right\}$, we have
\begin{align}
\frac{(-a)_{k+rp}(a+1)_{k+rp}}{(1)_{k+rp}^2}&\equiv \frac{(-a)_{r}(a+1)_{r}}{(1)_{r}^2}\cdot\frac{(-a)_{k}(a+1)_{k}}{(1)_{k}^2}\notag\\
&\hskip -10mm\times \left(1+2rpH_{\lfloor -pa \rfloor}-2rpH_k+rp\sum_{i=0}^{k-1}\left(\frac{1}{-a+i}+\frac{1}{a+1+i}\right)\right) \pmod{p^2}.\label{nc-1}
\end{align}
\end{lem}

\begin{lem}
For any even integer $n$, we have
\begin{align}
\sum_{k=0}^n{2k\choose k}^2{n+k\choose 2k}\left(-\frac{1}{4}\right)^k&=\frac{{n\choose n/2}^2}{4^n},\label{nc-2}\\
\sum_{k=0}^n{2k\choose k}^2{n+k\choose 2k}\left(-\frac{1}{4}\right)^kH_{k}&=\frac{{n\choose n/2}^2}{4^n}H_n,\label{nc-3}\\
\sum_{k=0}^n{2k\choose k}^2{n+k\choose 2k}\left(-\frac{1}{4}\right)^kH_{2k}&=\frac{{n\choose n/2}^2}{2\cdot4^n}H_n,\label{nc-4}\\
\sum_{k=0}^n{2k\choose k}^2{n+k\choose 2k}\left(-\frac{1}{4}\right)^k\sum_{i=0}^{k-1}\frac{1}{n+1+i}&=\frac{{n\choose n/2}^2}{4^n}\left(\frac{3}{2}H_n-H_{n/2}\right).\label{nc-5}
\end{align}
\end{lem}
See \cite[(2.9) \& (2.18)]{liujc-a07686-2016} and \cite[(16) \& (17)]{tauraso-a00729-2017}.

\begin{lem}
Let $p\ge 5$ be a prime. For $a\in\left\{-\frac{1}{2},-\frac{1}{3},-\frac{1}{4},-\frac{1}{6}\right\}$, we have
\begin{align}
H_{\lfloor-pa\rfloor}\equiv H_{\langle a\rangle_p}\pmod{p}.\label{nc-6}
\end{align}
\end{lem}
{\noindent\it Proof.}
For any prime $p\ge 5$, there exists $\varepsilon\in \{1,-1\}$ such that $p\equiv \varepsilon \pmod{2,3,4,6}$. We give the proof of \eqref{nc-6} for $a=-\frac{1}{3}$. The proofs of the other three cases run similarly.

If $p\equiv 1\pmod{3}$, then $\lfloor p/3\rfloor=(p-1)/3$ and $\langle -1/3\rangle_p=(p-1)/3$,
and so \eqref{nc-6} clearly holds for $a=-\frac{1}{3}$.

If $p\equiv -1\pmod{3}$, then $\lfloor p/3\rfloor=(p-2)/3$ and $\langle -1/3\rangle_p=(2p-1)/3$.
Applying the fact that $H_k=H_{p-1-k}\pmod{p}$ for $0\le k\le p-1$, we conclude that \eqref{nc-6} also holds for $a=-\frac{1}{3}$.
\qed

{\noindent\it Proof of Theorem \ref{nt-3}.}
We first prove that for any non-negative integer $r$,
\begin{align}
&\hskip -13mm\sum_{k=rp}^{(r+1)p-1}\frac{\left(\frac{1}{2}\right)_k(-a)_k(a+1)_k}{(1)_k^3}\notag\\
&\equiv (-1)^{\frac{p+1}{2}}\Gamma_p\left(-\frac{a}{2}\right)^2\Gamma_p\left(\frac{a+1}{2}\right)^2
\cdot\frac{\left(\frac{1}{2}\right)_r(-a)_r(a+1)_r}{(1)_r^3}\pmod{p^2}.\label{nc-7}
\end{align}
Noting that
\begin{align}
\frac{\left(\frac{1}{2}\right)_k(-a)_k(a+1)_k}{(1)_k^3}={2k\choose k}^2{a+k\choose 2k}\left(-\frac{1}{4}\right)^k,\label{nc-8}
\end{align}
and letting $k\to k+rp$ on the left-hand side of \eqref{nc-7}, we see that \eqref{nc-7} is equivalent to
\begin{align}
&\hskip-10mm\sum_{k=0}^{p-1}{2k+2rp\choose k+rp}^2{a+k+rp\choose 2k+2rp}\left(-\frac{1}{4}\right)^{k+rp}\notag\\
&\equiv (-1)^{\frac{p+1}{2}}\Gamma_p\left(-\frac{a}{2}\right)^2\Gamma_p\left(\frac{a+1}{2}\right)^2
{2r\choose r}^2{a+r\choose 2r}\left(-\frac{1}{4}\right)^r\pmod{p^2}.\label{nc-9}
\end{align}

By Lemma \ref{l1} and \ref{l2}, we have
\begin{align}
\left(\frac{1}{4}\right)^{k+rp}&=\left(\frac{1}{4}\right)^{k+r}\cdot \left(\frac{1}{4}\right)^{r(p-1)}\notag\\
&= \left(\frac{1}{4}\right)^{k+r}\left(1+pq_p\left(4^{-r}\right)\right)\notag\\
&\equiv \left(\frac{1}{4}\right)^{k+r}\left(1-2rpq_p\left(2\right)\right)\notag\\
&\equiv \left(\frac{1}{4}\right)^{k+r}\left(1+rpH_{\lfloor p/2\rfloor}\right)\pmod{p^2}.\label{nc-10}
\end{align}
Combining \eqref{nx} and \eqref{nc-10} gives
\begin{align}
&\hskip-10mm\left(\frac{1}{4}\right)^{k+rp}{2k+2rp\choose k+rp}\notag\\
&\equiv \left(\frac{1}{4}\right)^{k+r}{2r\choose r}{2k\choose k}(1+rp(2H_{2k}-2H_k+H_{\lfloor p/2\rfloor})\pmod{p^2}.\label{nc-11}
\end{align}
On the other hand, using \eqref{nc-6} and
\begin{align*}
\frac{(-a)_{k}(a+1)_{k}}{(1)_{k}^2}=(-1)^k{2k\choose k}{a+k\choose 2k},
\end{align*}
we can rewrite \eqref{nc-1} as
\begin{align}
&(-1)^{k+rp}{2k+2rp\choose k+rp}{a+k+rp\choose 2k+2rp}\notag\\
&\equiv (-1)^{r+k}{2r\choose r}{a+r\choose 2r}{2k\choose k}{a+k\choose 2k}\notag\\
&\times \left(1+2rpH_{\langle a\rangle_p}-2rpH_k+rp\sum_{i=0}^{k-1}\left(\frac{1}{-a+i}+\frac{1}{a+1+i}\right)\right) \pmod{p^2}.\label{nc-12}
\end{align}
Applying \eqref{nc-11} and \eqref{nc-12} to the left-hand side of \eqref{nc-9} yields
\begin{align}
&\text{LHS \eqref{nc-9}}\notag\\
&\equiv {2r\choose r}^2{a+r\choose 2r}\left(-\frac{1}{4}\right)^r
\sum_{k=0}^{p-1}{2k\choose k}^2{a+k\choose 2k}\left(-\frac{1}{4}\right)^k\notag\\
&\times \left(1+rp\left(2H_{2k}-4H_k+2H_{\langle a\rangle_p}+H_{\lfloor p/2\rfloor}+\sum_{i=0}^{k-1}\left(\frac{1}{-a+i}+\frac{1}{a+1+i}\right)\right)\right)
\pmod{p^2}.\label{nc-13}
\end{align}

In fact, the proof of \eqref{nc-9} follows from the congruence:
\begin{align}
&\sum_{k=0}^{p-1}{2k\choose k}^2{a+k\choose 2k}\left(-\frac{1}{4}\right)^k\notag\\
&\times\left(2H_{2k}-4H_k+2H_{\langle a\rangle_p}+H_{\lfloor p/2\rfloor}+\sum_{i=0}^{k-1}\left(\frac{1}{-a+i}+\frac{1}{a+1+i}\right)\right)
\equiv 0\pmod{p}.\label{nc-14}
\end{align}
Substituting \eqref{nc-14} into \eqref{nc-13}, and then using \eqref{nc-8} and \eqref{na-2}, we obtain
\begin{align*}
\text{LHS \eqref{nc-9}}
&\equiv {2r\choose r}^2{a+r\choose 2r}\left(-\frac{1}{4}\right)^r
\cdot\pFq{3}{2}{\frac{1}{2},-a,a+1}{1,1}{1}_{p}\\
&\equiv (-1)^{\frac{p+1}{2}}\Gamma_p\left(-\frac{a}{2}\right)^2\Gamma_p\left(\frac{a+1}{2}\right)^2
{2r\choose r}^2{a+r\choose 2r}\left(-\frac{1}{4}\right)^r\pmod{p^2},
\end{align*}
which is \eqref{nc-9}. So it suffices to prove \eqref{nc-14}.

Since $a\equiv \langle a\rangle_p\pmod{p}$, we have
\begin{align}
&\text{LHS \eqref{nc-14}}\notag\\
&\equiv \sum_{k=0}^{p-1}{2k\choose k}^2{\langle a\rangle_p+k\choose 2k}\left(-\frac{1}{4}\right)^k\notag\\
&\times\left(2H_{2k}-4H_k+2H_{\langle a\rangle_p}+H_{\lfloor p/2\rfloor}+\sum_{i=0}^{k-1}\left(\frac{1}{-\langle a\rangle_p+i}+\frac{1}{\langle a\rangle_p+1+i}\right)\right)
\pmod{p}.\label{nc-15}
\end{align}

By \eqref{nc-2}-\eqref{nc-5} and the fact that $\langle a\rangle_p$ is an even integer with $0\le \langle a\rangle_p\le p-1$, we obtain
\begin{align}
&\sum_{k=0}^{p-1}{2k\choose k}^2{\langle a\rangle_p+k\choose 2k}\left(-\frac{1}{4}\right)^k
\left(2H_{2k}-4H_k+2H_{\langle a\rangle_p}+H_{\lfloor p/2\rfloor}+\sum_{i=0}^{k-1}\frac{1}{\langle a\rangle_p+1+i}\right)\notag\\
&=\frac{{\langle a\rangle_p\choose \langle a\rangle_p/2}^2}{4^{\langle a\rangle_p}}
\left(H_{\lfloor p/2\rfloor}+\frac{1}{2}H_{\langle a\rangle_p}-H_{\langle a\rangle_p/2}\right).\label{nc-16}
\end{align}

Let $b=p-\langle a\rangle_p$. It is clear that $\langle a\rangle_p\equiv -b\pmod{p}$ and $b-1$ is an even integer with $0\le b-1\le p-1$. Using \eqref{nc-5} and the fact that ${-b+k\choose 2k}={b-1+k\choose 2k}$, we get
\begin{align}
&\sum_{k=0}^{p-1}{2k\choose k}^2{\langle a\rangle_p+k\choose 2k}\left(-\frac{1}{4}\right)^k
\sum_{i=0}^{k-1}\frac{1}{-\langle a\rangle_p+i}\notag\\
&\equiv \sum_{k=0}^{p-1}{2k\choose k}^2{-b+k\choose 2k}\left(-\frac{1}{4}\right)^k
\sum_{i=0}^{k-1}\frac{1}{b+i} \pmod{p}\notag\\
&=\sum_{k=0}^{p-1}{2k\choose k}^2{b-1+k\choose 2k}\left(-\frac{1}{4}\right)^k
\sum_{i=0}^{k-1}\frac{1}{b+i}\notag\\
&=\frac{{b-1\choose (b-1)/2}^2}{4^{b-1}}\left(\frac{3}{2}H_{b-1}-H_{(b-1)/2}\right).\label{nc-17}
\end{align}
By \eqref{nc-2}, we have
\begin{align}
\frac{{b-1\choose (b-1)/2}^2}{4^{b-1}}
&=\sum_{k=0}^{p-1}{2k\choose k}^2{b-1+k\choose 2k}\left(-\frac{1}{4}\right)^k\notag\\
&\equiv \sum_{k=0}^{p-1}{2k\choose k}^2{-\langle a\rangle_p-1+k\choose 2k}\left(-\frac{1}{4}\right)^k\pmod{p}\notag\\
&=\sum_{k=0}^{p-1}{2k\choose k}^2{\langle a\rangle_p+k\choose 2k}\left(-\frac{1}{4}\right)^k\notag\\
&=\frac{{\langle a\rangle_p\choose \langle a\rangle_p/2}^2}{4^{\langle a\rangle_p}}.\label{nc-18}
\end{align}

Substituting \eqref{nc-16}-\eqref{nc-18} into \eqref{nc-15} and then replacing $b$ by $b=p-\langle a\rangle_p$ gives
\begin{align}
\text{LHS \eqref{nc-14}}
&\equiv \frac{{\langle a\rangle_p\choose \langle a\rangle_p/2}^2}{4^{\langle a\rangle_p}}\notag\\
&\times\left(H_{\lfloor p/2\rfloor}+\frac{1}{2}H_{\langle a\rangle_p}-H_{\langle a\rangle_p/2}+\frac{3}{2}H_{p-1-\langle a\rangle_p}-H_{(p-1-\langle a\rangle_p)/2}\right)\pmod{p}.\label{nc-19}
\end{align}

Note that
\begin{align*}
H_{(p-1-\langle a\rangle_p)/2}&=H_{(p-1)/2}-\sum_{i=0}^{\langle a\rangle_p/2-1}\frac{1}{(p-1)/2-i}\notag\\
&\equiv H_{\lfloor p/2\rfloor}-\sum_{i=0}^{\langle a\rangle_p/2-1}\frac{1}{-1/2-i}\pmod{p}\notag\\
& = H_{\lfloor p/2\rfloor}+2H_{\langle a\rangle_p}-H_{\langle a\rangle_p/2}.
\end{align*}
Substituting the above congruence into \eqref{nc-19} yields
\begin{align*}
\text{LHS \eqref{nc-14}}\equiv \frac{3}{2}\cdot\frac{{\langle a\rangle_p\choose \langle a\rangle_p/2}^2}{4^{\langle a\rangle_p}} \left(H_{p-1-\langle a\rangle_p}-H_{\langle a\rangle_p}\right)\equiv 0 \pmod{p},
\end{align*}
where we have utilized the fact that $H_{p-1-k}\equiv H_k\pmod{p}$ for $0\le k \le p-1$.
This concludes the proof of \eqref{nc-14}.

Taking the sum over $r$ from $0$ to $n-1$ on both sides of \eqref{nc-7} gives
\begin{align*}
&\hskip -5mm\pFq{3}{2}{\frac{1}{2},-a,a+1}{1,1}{1}_{np}\\
&\equiv (-1)^{\frac{p+1}{2}}\Gamma_p\left(-\frac{a}{2}\right)^2\Gamma_p\left(\frac{a+1}{2}\right)^2
\pFq{3}{2}{\frac{1}{2},-a,a+1}{1,1}{1}_{n} \pmod{p^2}.
\end{align*}
This completes the proof of Theorem \ref{nt-3}.\qed

\vskip 5mm \noindent{\bf Acknowledgments.} The author would like to
thank the anonymous referee for careful reading of this manuscript and helpful comments.

\end{document}